\documentclass[twocolumn]{article}

\usepackage{graphicx}

\newcommand{\be}{\begin{equation}}
\newcommand{\ee}{\end{equation}}

\title{\bf Asymptotic Analysis of the Paradox \\
in Log-Stretch Dip Moveout}

\author{Xin-She Yang \\
{\small Department of Applied Mathematics,
University of Leeds,
Leeds LS2 9JT, UK.   }
\and
Binzhong Zhou \\
{\small CSIRO Division of Exploration and Mining, 
P.O.Box 883, Kenmore, QLD 4069, AUSTRALIA. }}

\date{}

\begin{document}
\maketitle

\begin{abstract}

There exists a paradox in dip moveout (DMO) in seismic  data processing.
The paradox is why  Notfors and Godfrey's approximate time
log-stretched DMO can produce better impulse responses than the
full log DMO, and why Hale's f-k DMO is correct although
it was based on two inaccurate assumptions for the midpoint
repositioning and the DMO time relationship? Based on the asymptotic
analysis of the DMO algorithms,  we find that any form of correctly
formulated DMO must handle both space and time coordinates properly in order
to deal with all dips accurately. The surprising improvement of
Notfors and Godfrey's log DMO on  Bale and Jakubowicz's full log
DMO was due to the equivalent midpoint repositioning by transforming
the time-related phase shift to the space-related phase shift.
The explanation of why Hale's f-k DMO is correct although it was
based on two inaccurate assumptions is that the two approximations
exactly cancel each other in the f-k domain to give the correct
final result. \\

\end{abstract}

\noindent {\bf Citation detail:} X. S. Yang and B. Zhou,
Asymptotic analysis of the paradox in log-stretch dip moveout,
{\it Geophys. Res. Lett.}, {\bf 27}(3), 441-444 (2000).

\section*{INTRODUCTION}

Dip Moveout (DMO) technique has been widely used in
seismic data processing over
the past decade and many different algorithms have been developed [{\it
Hale, 1984; Bale and Jakubowicz, 1987; Notfors and Godfrey, 1987;
Liner, 1990; Gardner, 1991; Black et al., 1993; Zhou et al., 1996}].
Hale's DMO is accurate for all reflector dips and has
become an industrial standard, but it is computational
expensive and temporally non-stationary. In all the methods, a
logarithmic time-stretching technique [{\it Bolondi et al., 1982}]
is widely used, probably due to its computational efficiency and easy
implementation. The impulse response produced by Bale and
Jakubowicz's  full log-stretch operator (hereafter referred to as
Bale's full log DMO) in the frequency-wavenumber (f-k) domain is inaccurate
although no approximations were made in the mathematical
formulations.  The impulse response is surprisingly improved by
Notfors and Godfrey's approximate log-stretch scheme (hereafter referred
to as Notfors' log DMO). Bale's inaccurate DMO impulse responses reveal that
Hale's DMO derivation had some  inappropriate assumptions.
To improve the DMO impulse responses,  an exact log-stretch DMO
was derived by transforming Hale's time log-stretch impulse
response into the Fourier domain [{\it Liner, 1990}] .
{\it Black et al.} [1993] provided another
derivation of the more reliable f-k DMO by treating
the reflection-point smear from dipping reflectors correctly.
Based on Black's DMO relationship, {\it Zhou et al.} [1995, 1996]
presented a log-stretch f-k DMO which
is  similar to Gardner's DMO [{\it Gardner, 1991}] but  from a
derivation with more direct physical insight into the DMO process.
All these   DMO schemes  were developed to try to generate
better impulse responses and to improve the computational speed, but
have  different corrections of  Hale's
subtle flaw and improving Bale's inaccurate impulse responses.

Despite extensive studies and the routine use of the DMO algorithms,
there still exists a paradox which has not yet been explained
satisfactorily. The paradox is why Notfors' approximate
log DMO can produce better impulse responses than Bale's full log
DMO, and why Hale's f-k DMO
is correct although it was based on two inaccurate assumptions.
Liner pointed out [{\it Liner, 1990}] that Bale and Jakubowicz's
DMO derivation implicitly assumes that the Fourier transform frequency
in the log-stretch domain is time-independent. Black et al.
attributed Hale's subtle flaw to the lack of midpoint change
and the consequently
improper treatment of the reflection-point smear from dipping reflectors.
Zhou et al. suggested that the approximations from the non-repositioned
midpoint  and from the incorrect time
relationship before and after DMO  may counteract
each other in the f-k domain to make the final result correct, but they
may not cancel each other in the log-stretch f-k domain. In this paper,
we analyse these viewpoints and compare
DMO algorithms to find out what is really responsible for
the inaccurate impulse responses of Bale's full log DMO scheme.
For convenience in the following discussion, only common-offset DMO
forms will be analysed although  similar methodology can be also
applied to common-shot DMO analysis.

\section*{ASYMPOTIC ANALYSIS OF PHASE SHIFT FUNCTIONS}

We assume a constant-velocity medium
and follow the notations used earlier [{\it Hale, 1984;
Notfors and Godfrey, 1987; Black et al., 1993; Zhou et al., 1996}].
The heuristic DMO mapping can be expressed by
\be
P_{0} (t_{0}, x_{0}, h)=P_{n} (t_{n}, x_{n}, h),    \label{DMO:MAPPING}
\ee
where $P_{n} (t_{n}, x_{n},h)$ is a normal-moveout (NMO)
corrected input section and $P_{0} (t_{0}, x_{0}, h)$ is a zero-offset
DMO output section. $h$ is the half offset,
$t_{n}$ and $x_{n}$ are time and space  in the NMO-corrected domain
respectively, and  $t_{0}$ and $x_{0}$ are the time and space in
DMO-corrected domain, respectively. The time log-stretch transform pair
[{\it Bolondi et al., 1982; Bale \& Jakubowicz 1987}]
\be
\tau = {\rm ln} (t/t_{c}) \,\,\,\,\,\, {\rm and}\,\,\,\,\,\,
t=t_{c} e^{\tau},    \label{TAU:LOGT}
\ee
transforms the time coordinates from $(t_{0}, t_{n})$ to the log-stretch
variables $(\tau_{0}, \tau_{n})$. $t_{c}$ is the minimum cutoff time
introduced to avoid the singularity of the logarithm at zero.

The DMO processing algorithm can be divided into the following five steps
[{\it Liner, 1990; Gardner, 1991; Hale, 1991}]: 1) transform NMO
corrected input
data $P_{n} (t_{n}, x_{n},h)$ to $P_{n} (\tau_{n}, x_{n},h)$ by the time
log stretch relation (\ref{TAU:LOGT}), 2) 2-D FFT of $P_{n} (\tau_{n}, x_{n},
h)$ over $\tau_{n}$ and $x_{n}$ to $P_{n} (\Omega, k, h)$, 3) use the
phase factor for filtering in the $(\Omega, k)$ domain to get
$P_{0} (\Omega, k, h)$, 4) inverse 2-D FFT of $P_{0} (\Omega, k, h)$ over
$\Omega$ and $k$ to $P_{0} (\tau_{0}, x_{0}, h)$, and 5) inverse
log stretch of $P_{0} (\tau_{0}, x_{0}, h)$ by the inverse relation
of equation (\ref{TAU:LOGT}) to get the  output $P_{0} (t_{0}, x_{0}, h)$.
Of all the available log stretch DMO algorithms, the same transform
relation (\ref{TAU:LOGT})  is used in steps 1) and 5), and the
two-dimensional FFT transform pair is also used in steps 2) and 4). The
main difference between algorithms is thus the multiplying phase factor in
step 3).  The following discussion will therefore focus on
the analysis of the geometry-determining phase factor in different
log-stretch DMO algorithms.
As we will see, all the phase shift functions depend mainly on the  variable
\be
\xi=\frac{h k}{\Omega}.
\ee
Mathematically speaking, $\xi$ can
vary from $0$ (zero dip at the centre of the DMO ellipse) to
$\infty$ (very shallow, steeply dipping).

\subsection*{Bale's Full Log DMO}

From Hale's DMO formulation, {\it Bale and Jakubowicz} [1987] derived the
 full log stretch DMO
\be
P_{0} (\Omega, k, h)=e^{-i \frac{1}{2} \Omega {\rm ln}(1-\xi^{2}) }
P_{n} (\Omega, k, h),  \label{BALE:DMO}
\ee
However, this algorithm does not produce good
impulse responses. {\it Liner } [1990] attributed this feature to the fact that
Bale and Jakubowicz's derivation implicitly assumes that $\Omega$ is
{\em not} time-dependent. In fact, Bale's phase shift function
\be
\Phi_{F}=-\frac{1}{2} \Omega {\rm ln}(1-\xi^{2}), \label{FULL}
\ee
is only physically meaningful when $\xi < 1$ (near zero dips). It becomes
complex if $\xi \ge 1$ (shallow, steep dips) which is responsible for
the {\em inverted Gaussian} shape.

\subsection*{Notfors' Log DMO}

An  approximation to (\ref{BALE:DMO}) was derived by {\it Notfors and
Godfrey} [1987]
\be
P_{0} (\Omega, k, h)=e^{i \Omega ( \sqrt{1+\xi^{2}}-1 )}
P_{n} (\Omega, k, h),
\ee
which is based on two approximations:
(1) the independence of $\Omega$ and $\tau_{n}$,
and (2) $(h k/\Omega)^{2} << 1$. With these approximations,
Bale's phase shift function becomes Notfors' phase shift function
\be
\Phi_{N}=\Omega ( \sqrt{1+\xi^{2}}-1 ), \label{NOTFORS}
\ee
which is meaningful for any value of $\xi$ (all reflector dips).
This derivation implicitly extends the definition range of $\xi$ although
the approximations are based on $\xi << 1$.
It is this second approximation that makes Notfors' log DMO produce better
impulse responses than Bale's full log DMO (\ref{BALE:DMO}).

\subsection*{Liner's Log DMO}

The exact log stretch DMO derived by {\it Liner} [1990] is
\be
P_{0} (\Omega, k, h)=\frac{e^{i [\Omega \Delta_{s}-k y_{s}]}}{(1+
\beta_{s}^{2})^{1/2}} P_{n} (\Omega , k, h),  \label{LINER:DMO}
\ee
where
\be
\beta_{s}=y_{s}/h,
\ee
\be
\Delta_{s}=\frac{1}{2} {\rm ln} (1-\beta_{s}),
\ee
and the stationary point is
\be
y_{s}=\frac{h}{2 \xi} (1-\sqrt{1+4 \xi^{2}}).  \label{LINER:YS}
\ee
Liner's derivation introduces  log stretch variables into
Hale's $(t,x)$ elliptical response rather than into the NMO equation,
as was done earlier by {\it Bale and Jakubowicz} [1987].  The formula
(\ref{LINER:DMO}) does yield the correct impulse response geometry
[{\it Liner, 1990; Zhou et al., 1996}].

\subsection*{Zhou's Log DMO}

{\it Zhou et al.} [1996] presented an accurate log stretch DMO
\be
P_{0}(\Omega,k,h)=e^{i \frac{1}{2} \Omega [ \sqrt{1+4 \xi^{2}}
-1 -{\rm ln}[\frac{1}{2} (\sqrt{1+4 \xi^{2}}+1) ] } \\
P_{n}(\Omega, k, h),   \label{ZHOU:DMO}
\ee
which is equivalent to Gardner's results except for a sign difference
on the phase term due to a different definition of the 2-D Fourier transform
in Gardner's formulation. {\it Zhou et al.} [1996] derived their result
based on the DMO relationships given by {\it Black et al.} [1993].
It is easy to check that the phase shift in equation
(\ref{ZHOU:DMO}) is identical to that in equation (\ref{LINER:DMO}).
Therefore, Zhou's log DMO and Liner's log DMO will produce the same
impulse response geometry. The only difference is that Zhou's DMO
will yield slightly larger
amplitudes than Liner's DMO because (\ref{LINER:DMO}) has an amplitude
factor $1/(1+\beta_{s}^{2})^{1/2} \le 1$. For convenience in the following
analysis, we will refer to Liner's log DMO and Zhou's log DMO as
the exact log DMO, and will use the simpler phase shift expression
\be
\Phi_{E}=\frac{1}{2} \Omega [ \sqrt{1+4 \xi^{2}}
-1 -{\rm ln}[\frac{1}{2} (\sqrt{1+4 \xi^{2}}+1) ], \label{EXACT}
\ee
which is accurate for all reflector dips.

All the phase shift functions (\ref{FULL}),
(\ref{NOTFORS}) and (\ref{EXACT}) are monotonically increasing, therefore,
we  discuss only two asymptotic cases.

\subsection*{Case 1: $\xi << 1$}

For $\xi \rightarrow 0$, equations (\ref{FULL}), (\ref{NOTFORS}),
(\ref{EXACT}) all become
\be
\Phi_{F,N,E} \approx \frac{1}{2} \Omega \xi^{2},
\ee
which means that all these log-stretch DMO algorithms are asymptotically
equivalent near zero dips ($\xi \rightarrow 0$) and will yield virtually
the same impulse response geometry near the center of the DMO ellipse.

\subsection*{Case 2: $\xi >> 1$}

As the reflector dips become shallow and steep (for $\xi >> 1$), expression
(\ref{FULL}) will be invalid, but (\ref{NOTFORS}) and (\ref{EXACT})
can be asymptotically approximated as
\be
\Phi_{N} \approx \Omega \xi, \label{EXACT:N}
\ee
and
\be
\Phi_{E} \approx \Omega \xi (1-\frac{\ln \xi}{2 \xi}), \label{EXACT:E}
\ee
respectively. Using $\frac{\ln \xi}{\xi} \rightarrow 0$ as
$\xi \rightarrow \infty$, we can write (\ref{EXACT:N})
and (\ref{EXACT:E}) as a single asymptotic form
\be
\Phi_{N,E} \sim \Omega \xi,
\ee
which means that Notfors' log-stretch DMO and the exact log-stretch are
{\em asymptotically} equivalent but they are {\em not} approximately
equal due to the term $O(\frac{\ln \xi}{\xi})$ in (\ref{EXACT:E}).
This implies that the exact log DMO is able to deal with all
reflector dips correctly, but Notfors' log DMO mishandles the shallow,
steep dips although it is an improvement on Bale's DMO scheme.

\section*{THE EXPLANATION OF THE PARADOX}

Based on the above analysis, we can explain  the paradox of
why Notfors approximate log DMO improved Bale's full log DMO, and find
out the real flaw in Hale's formulation. The f-k integral for common-offset
in {\it Hale's} DMO [1984] is
\[ P_{0} (\omega,k,h) = \int \! \int P_{0} (t_{0}, x_{0}, h)
e^{i (\omega t_{0} - k x_{0}) } dx_{n} dt_{n} \]
\be
= \int \!\int \frac{1}{A} e^{i (\omega t_{n} A - k x_{n}) }
P_{n} (t_{n}, x_{n}, h) dx_{n} dt_{n},         \label{HALE:DMO}
\ee
where
\be
A= \sqrt{1+(\frac{h}{t_{n}} \frac{dt_{0}}{dx_{0}})^{2}}
= \sqrt{1+(\frac{h k}{t_{n} \omega} )^{2}}.   \label{DMO:A}
\ee
and the DMO mapping (\ref{DMO:MAPPING}) and the following relations of
time and space coordinate are used
\be
t_{0}= A t_{n} \,\,\,\,\,\,\,\,\,\, {\rm and}
\,\,\,\,\,\,\,\,\,\,x_{0}=x_{n},  \label{HALE:TX}
\ee
which are not strictly correct and make the DMO operator amplitude-unpreserved.
The similar time and space relationships were used by {\it Bale and
Jakubowicz} [1987] in their full log DMO derivation. {\it Black et al.}
[1993] derived a new amplitude preserving DMO integral
\[
P_{0} (\omega, k, h) = \int \! \int P_{0} (t_{0}, x_{0}, h) e^{i (\omega t_{0} - k x_{0}) } dx_{n} dt_{n} \]
\be
= \int \! \int  \frac{2 A^{2}-1}{A^{3}}
e^{i (\omega t_{n} A - k x_{n}) } P_{n} (t_{n}, x_{n}, h)
dx_{n} dt_{n}.     \label{BLACK:DMO}
\ee
by substituting (\ref{DMO:MAPPING}) and basing  on Black's correct
relationship of time and space coordinates before and after DMO processing
\be
t_{0}=\frac{t_{n}}{A} \,\,\,\,\,\,\,\,\,\, {\rm and}
\,\,\,\,\,\,\,\,\,\,x_{0}=x_{n}-\frac{1}{A} \frac{h^{2}}{t_{n}}
\frac{k}{\omega},  \label{BLACK:TX}
\ee
where the definition of $A$ is the same as in (\ref{DMO:A}). From
equation (\ref{BLACK:TX}), we can easily derive the DMO impulse
ellipse (Zhou et al., 1995) while the relationships in equation
(\ref{HALE:TX})give an impulse point.  As $A \ge 1$,
Black's correct relationship (\ref{BLACK:TX}) implies the inequality
\be
t_{0} \le t_{n}.  \label{BLACK:T0TN}
\ee
It is clearly seen that phase shifts in Black's DMO
(\ref{BLACK:DMO}) and in Hale's DMO (\ref{HALE:DMO}) are identical
[{\it Black et al., 1993}] due to the identical phase factor
$\exp[ i (\omega t_n A-k x_n)]$.  This implies that the two
approximations resulting from the inaccurate
relationships (\ref{HALE:TX}) exactly cancel each other in the f-k domain,
and thus make the final result correct. However, the counteraction of the
two approximations does not give exact cancellation in Bale's
log-stretch DMO algorithm because only the {\em time} coordinate is
stretched and the {\em space} coordinate remains unchanged. Therefore,
we expect that the phase shift in Bale's log DMO is
essentially due to the log-stretch traveltime relationship
\be
\tau_{0}-\tau_{n}=-\frac{1}{2} {\rm ln} (1-\xi^{2}). \label{BALE:T0TN}
\ee
This implies the incorrect inequalities
\be
\tau_{0} \ge \tau_{n}\,\,\,\,\,\,\,\,\,\,{\rm and}
\,\,\,\,\,\,\,\,\,\, t_{0} \ge t_{n},  \label{INEQU}
\ee
which contradict the correct relation (\ref{BLACK:T0TN}). {\it Notfors
and Godfrey} [1987] used the approximate version
of (\ref{BALE:T0TN})
\be
\tau_{0}-\tau_{n}= \sqrt{1+\xi^{2}}-1, \label{NOTFORS:T0TN}
\ee
which yields the same incorrect inequality as in (\ref{INEQU}).

From the general forms of log-stretch DMO algorithms,
we know that the phase shift function
\be
\Phi= k (x_{n}-x_{0}) -\Omega (\tau_{n} -\tau_{0}),
\ee
consists of two parts: a space-related phase shift $k (x_{n}-x_{0})$
and a time-related phase shift $\Omega (\tau_{n} -\tau_{0})$. The asymptotic
analysis in the previous section reveals that the space-related
phase shift is always greater than the time-related phase shift.
Comparing the phase shift functions (\ref{NOTFORS}) and (\ref{EXACT}),
we see that Notfors approximate log DMO virtually transformed
the time-related phase shift into its space-related counterpart,
and thus equivalently repositioned the midpoint
\be
x_{0}=x_{n}+\frac{h}{\xi} (1-\sqrt{1+\xi^{2}}),
\ee
where the second term on the right hand side is similar to the
expression (\ref{LINER:YS}) of the stationary point $y_{s}$ in
Liner's exact log DMO. It is this
equivalent midpoint repositioning that makes Notfors' approximate
log DMO overcome Hale's subtle flaw and subsequently produces a
better match to the geometry of the impulse response DMO ellipse,
but the Notfors' wider geometry of the DMO ellipse shows that
the shallow steep dips are still improperly treated.

\section*{CONCLUSIONS}

The analysis  of different log-stretch DMO algorithms
has shown that any form of correctly formulated DMO must handle
both space and time coordinates properly in order to deal with
all dips accurately. The log-stretch DMO algorithms presented by
{\it Liner} [1990] and
{\it Zhou et al.} [1996] are able to treat all reflector dips correctly.
Bale and Jakubowicz's full log DMO is only
approximately correct for the near zero dips, and it becomes
singular for shallow, steep dips. Notfors and Godfrey's approximate
log DMO is a great improvement in handling all dips, but
its impulse response still departs from the exact DMO ellipse,
especially for shallow steep dips. The surprising improvement
of Notfors and Godfrey's log DMO on Bale and Jakubowicz's full
log DMO was mainly due to the equivalent midpoint repositioning
by equivalently transforming the time-related phase shift to the
space-related phase shift. The explanation of why Hale's f-k DMO
is correct although it was based on two inaccurate assumptions
is that the two approximations exactly cancel each other in the
f-k domain to give the correct  final result.

{\bf Acknowledgements.}  The authors wish to thank the anonymous referees
for their very helpful comments. This research was supported by CSIRO
Exploration and Mining in Australia and the Cooperative Research Centre
for Mining Technology and Equipment.

\section*{REFERENCES}

\begin{description}
\item Bale, R., and Jakubowicz, H.,  Post-stack prestack migration,
57th Ann. Mtg., {\it Soc. Exp. Geophys.}, Expanded Abstract, 714-717, 1987.

\item Black, J. L., Schleicher, K. L., and Zhang, L., True-amplitude imaging
and dip moveout: {\it Geophysics}, {\it 58}, 47-66, 1993.

\item Bolondi, G., Loinger, E., and Rocca, F.,  Offset continuation
of seismic sections, {\it Geophys. Prosp.}, {\it 30}, 813-828, 1982.

\item Gardner, G. H. F., Interpolation, crossline migration and inline
depth migration of 3-D marine surveys, {\it The SAL Annual Progress Review},
{\it 24}, 57-69, 1991.

\item Hale, D.,  Dip moveout by Fourier transform, {\it Geophysics},
{\it 49}, 741-757, 1984.

\item Liner, C. L.,  General theory and comparative anatomy of dip moveout,
{\it Geophysics}, {\it 55}, 595-607, 1990.

\item Notfors, C. D., and Godfrey, R, J., Dip moveout in the
frequency-wavenumber domain, {\it Geophysics}, {\it 52}, 1718-1721, 1987.

\item Zhou, B., Mason, I. M., and Greenhalgh, S. A., An accurate
formulation of log-stretch dip moveout in the frequency-wavenumber
domain, {\it Geophysics}, {\it 61},  815-821, 1996.

\item Zhou, B., Mason, I. M., and Greenhalgh, S. A., Accurate and efficient
shot-gather dip moveout processing in the log-stretch domain, {\it
Geophysical Prospecting}, {\it 43}, 863-978, 1995.

\end{description}

\end{document}